\documentclass[a4paper,reqno,11pt]{amsart}
\usepackage[margin=1in]{geometry}
\usepackage[T1]{fontenc}
\usepackage{lmodern}
\usepackage{amsmath,amssymb,amsthm,mathtools}
\usepackage{microtype}
\usepackage[foot]{amsaddr}
\usepackage{xcolor}
\usepackage[colorlinks=true,linkcolor=blue!45!black,
  citecolor=green!35!black,urlcolor=blue!60!black]{hyperref}
\usepackage[nameinlink,noabbrev]{cleveref}

\newtheorem{theorem}{Theorem}[section]
\newtheorem{lemma}[theorem]{Lemma}
\newtheorem{proposition}[theorem]{Proposition}
\newtheorem{corollary}[theorem]{Corollary}
\newtheorem{conjecture}[theorem]{Conjecture}
\theoremstyle{remark}

\newcommand{\Z}{\mathbb Z}
\newcommand{\F}{\mathbb F}
\newcommand{\cP}{\mathcal P}
\newcommand{\1}{\mathbf 1}
\allowdisplaybreaks[1]
\title{On ill-distributed $B_3$ and $B_4$ sets}
\author{Andre Hisatsuga}
\author{Rafael Miyazaki}
\author{Cosmin Pohoata}
\address{Department of Mathematics, Emory University, Atlanta, GA 30322, USA}
\email{\{andre.hisatsuga|rafael.kazuhiro.miyazaki|cosmin.pohoata\}@emory.edu}

\date{}

\begin{document}
\begin{abstract}
We show that, for $k\in\{3,4\}$, a $B_k$-set
$A\subseteq[N]$ has size $o(N^{1/k})$ if it omits a fixed positive
proportion of the residue classes modulo every prime in a fixed family
whose reciprocals have divergent sum.
\end{abstract}
\maketitle

\section{Introduction}

For an integer $k\ge2$, a set $A$ of integers is a \emph{$B_k$-set} if
every equality
\[
 a_1+\cdots+a_k=b_1+\cdots+b_k,
 \qquad a_i,b_i\in A,
\]
forces the two lists of summands to agree as multisets. Repetitions are
allowed. The case $k=2$ is the more familiar notion of a \emph{Sidon set}.
Writing $[N]=\{1,\ldots,N\}$ and
$kA=\{a_1+\cdots+a_k:a_i\in A\}$, a $B_k$-set $A\subseteq[N]$
satisfies
\begin{equation}\label{eq:counting}
 \frac{|A|^k}{k!}\le \binom{|A|+k-1}{k}=|kA|\le kN.
\end{equation}
Thus $|A|=O_k(N^{1/k})$. For a subset of
$\{1^k,2^k,\ldots,n^k\}$, this gives only $O_k(n)$; obtaining $o(n)$
requires using the arithmetic of the ambient set.

For a prime $p$, let
\[
 A_p=\{a\bmod p:a\in A\}\subseteq\Z/p\Z.
\]
Sets occupying relatively few residue classes are often called
\emph{ill-distributed}. A standard consequence of Gallagher's larger
sieve~\cite{Gallagher} is that
\[
 |A_p|\le\alpha p+C\quad\text{for every prime }p
 \qquad\Longrightarrow\qquad
 |A|\ll_{\alpha,C}N^\alpha
\]
for fixed $0<\alpha<1$ and $C\ge0$. At $\alpha=1/k$, this is the
same scale as the elementary $B_k$ bound~\eqref{eq:counting}. These bounds arise from different restrictions: additive uniqueness and scarcity of occupied residue classes. It turns out that imposing both can force a further saving.

In a recent paper, Croot, Mao, Pohoata, Sheffer, and Yip~\cite[Theorem~1.1]{CMPSY}
proved this quantitatively for Sidon sets: if
$A\subseteq[N]$ is Sidon and $|A_p|\le\alpha p$ for every prime,
then, for some $c_\alpha>0$,
\begin{equation} \label{Sidon}
 |A|\ll_\alpha N^{1/2}
 \exp\!\left(-c_\alpha\frac{\log N}{\log\log N}\right).
\end{equation}
The case $\alpha=1/2$, up to a harmless bounded error in $|A_p|$,
connects this result to the classical problem of Alon and
Erd\H{o}s~\cite{AE} on the largest Sidon subset of the squares.
Indeed, if $A\subseteq\{1^2,2^2,\ldots,n^2\}$, then
$|A_p|\le(p+1)/2$ for every odd prime $p$. Taking $N=n^2$ gives
\[
 |A|\ll n\exp\!\left(-c\frac{\log n}{\log\log n}\right)
\]
for every such Sidon set $A$ and some absolute constant $c>0$;
see~\cite[Corollary~1.2]{CMPSY}. 

The simple sieve theoretic ideas behind this result, sometimes in combination with entropic methods or with the theory of towers of number fields with bounded root discriminant~\cite{Martinet,HM,HMR} (recently brought to light through the breakthrough OpenAI
unit-distance construction~\cite{OpenAI,UnitRemarks}), have led to several new and unexpected results in combinatorics, number theory, and
discrete geometry. See, for example, \cite{CMPSY, PohoataER,LPZ,CMPY}.

In the same first paper \cite{CMPSY}, the authors also proved that every
$B_3$-set $A\subseteq\{1^3,2^3,\ldots,n^3\}$ satisfies
\[
 |A|\ll n\exp\!\left(-c_3\frac{\sqrt{\log n}}{\log\log n}\right),
\]
while every $B_4$-set $A\subseteq\{1^4,2^4,\ldots,n^4\}$ satisfies
\[
 |A|\ll\frac{n}{(\log\log n)^{c_4}},
\]
for absolute constants $c_3,c_4>0$. These results are considerably more
involved than the argument for Sidon sets, relying on a rather different entropic sieve method, one which also leverages further arithmetic information specific to the cubes and fourth powers, respectively (beyond ill-distribution). We refer to \cite{CMPSY} for more details.

The main contribution of this paper is a completely elementary method giving $|A|=o(N^{1/3})$
for $B_3$-sets and $|A|=o(N^{1/4})$ for $B_4$-sets in $[N]$ under an ill-distribution hypothesis alone. More precisely, in
addition to the $B_k$ property, it suffices that $A$ misses
\emph{any} fixed positive proportion of residue classes modulo every
prime in a fixed family with divergent reciprocal sum.

\begin{samepage}
\begin{theorem}\label{thm:main}
Let $k\in\{3,4\}$, let $0<\alpha<1$, and let $\cP$ be a fixed family
of primes such that
\[
 \sum_{p\in\cP}\frac1p=\infty.
\]
If $A\subseteq[N]$ is a $B_k$-set and $|A_p|\le\alpha p$ for every
$p\in\cP$, then
\[
 |A|=o_{\alpha,\cP}(N^{1/k})\qquad(N\to\infty).
\]
\end{theorem}
\end{samepage}

The parameters $\alpha$ and $\cP$ are fixed independently of $N$ and
$A$. Removing finitely many primes from $\cP$ makes no difference.
Consequently, the hypothesis may also be replaced by
$|A_p|\le\alpha p+C$ for a fixed constant $C$: choose
$\alpha<\alpha'<1$ and discard the finitely many primes for which
$\alpha p+C>\alpha'p$.

It is also perhaps worth mentioning that for Sidon sets $A\subseteq[N]$, the qualitative conclusion $|A|=o(N^{1/2})$ already follows from the proof of \cite[Theorem~1.1]{CMPSY}, if we only require $|A_p|\le\alpha p$ to hold on any fixed infinite family $\cP$ of primes, with fixed $0<\alpha<1$ (in particular, no divergence assumption is even needed). 

As already alluded to, Theorem \ref{thm:main} recovers the qualitative conclusions from \cite{CMPSY} for cubes and
fourth powers.

\begin{corollary}\label{cor:powers}
For each $k\in\{3,4\}$, every $B_k$-set
$A\subseteq\{1^k,2^k,\ldots,n^k\}$ satisfies $|A|=o(n)$.
\end{corollary}

The proof of Theorem \ref{thm:main} will in some sense extend the original (simple) proof of \eqref{Sidon}. The main idea will be to consider the \emph{modular additive energy}
\[
 E_Q(A)=\#\{(a,b,c,d)\in A^4:a-b\equiv c-d\pmod Q\},
\]
for a suitable squarefree integer $Q$ made up of primes from $\mathcal{P}$. The easy bound $E_Q(A)\ge |A|^4/Q$ follows from a direct application of the Cauchy-Schwarz inequality. Using a much more delicate Cauchy-Schwarz argument, however, we can show that
\begin{equation}\label{eq:intro-energy}
 E_Q(A)\ge\frac{|A|^4}{Q}
 \prod_{p\mid Q}\left(1+
 \frac{(p/|A_p|-1)^2}{p-1}\right)
\end{equation}
holds for every squarefree integer $Q$. The ill-distribution assumption on $A$ then gives us a factor of $1+\Omega_\alpha(1/p)$ at each prime $p|Q$, and then the product can be made arbitrarily large. Elementary $B_3$ and $B_4$ counting bounds give the opposing
upper estimate. For the qualitative conclusion, one can fix $Q$,
let $N$ tend to infinity, and then enlarge $Q$.

\section{A lower bound for $E_Q(A)$}

\begin{proposition}\label{prop:lower}
Let $A$ be a nonempty finite set of integers and let $Q$ be a
squarefree positive integer. Then~\eqref{eq:intro-energy} holds.
\end{proposition}

The idea is that a small residue image forces many differences to be
divisible by each prime. The weights below combine these gains without
any assumption that the residue coordinates of $A$ are independent.

\begin{proof}
For each residue class $t\bmod Q$, write
\[
 r(t)=\#\{(a,b)\in A^2:a-b\equiv t\pmod Q\}.
\]
Since $\sum_{t\bmod Q}r(t)=|A|^2$, the usual Cauchy--Schwarz
inequality gives
\[
 E_Q(A)=\sum_{t\bmod Q}r(t)^2
 \ge\frac{\bigl(\sum_{t\bmod Q}r(t)\bigr)^2}{Q}
 =\frac{|A|^4}{Q}.
\]
To exploit the missing residue classes, we replace the constant weight
by one that favors differences divisible by primes dividing $Q$.
For any nonzero real function $w$ on $\Z/Q\Z$, Cauchy--Schwarz gives
\begin{equation}\label{eq:weighted-cs}
 E_Q(A)\ge
 \frac{\bigl(\sum_{t\bmod Q}r(t)w(t)\bigr)^2}
      {\sum_{t\bmod Q}w(t)^2}.
\end{equation}
For each prime $p\mid Q$ and each residue class $t\bmod Q$, define
\[
 w_p(t)=|A_p|-1+(p-|A_p|)\1_{p\mid t}
 =\begin{cases}p-1,&p\mid t,\\|A_p|-1,&p\nmid t,\end{cases}
 \qquad
 w(t)=\prod_{p\mid Q}w_p(t).
\]
Here $\1_{p\mid t}$ equals $1$ if $t\equiv0\pmod p$ and $0$
otherwise; this is well-defined for $t\bmod Q$ since $p\mid Q$.
Thus $w_p$ assigns an extra weight of $p-|A_p|$ to differences
divisible by $p$, and $w$ combines these local weights.
We now estimate the two sums in~\eqref{eq:weighted-cs}.

We first bound $\sum_t r(t)w(t)$ from below.
For every $d\mid Q$, the Chinese remainder theorem shows that $A$
occupies at most $\prod_{p\mid d}|A_p|$ residue classes modulo $d$.
Cauchy--Schwarz applied to their occupancies gives
\begin{equation}\label{eq:divisor-collisions}
 \sum_{\substack{t\bmod Q\\d\mid t}}r(t)
 =\#\{(a,b)\in A^2:a\equiv b\pmod d\}
 \ge\frac{|A|^2}{\prod_{p\mid d}|A_p|}.
\end{equation}
All coefficients in the expansion of $w$ as a linear combination of
the functions $\1_{d\mid t}$ are nonnegative. We may therefore
apply~\eqref{eq:divisor-collisions} term by term to obtain
\begin{align}
 \sum_{t\bmod Q}r(t)w(t)
 &\ge |A|^2\sum_{d\mid Q}
   \prod_{p\mid Q/d}(|A_p|-1)
   \prod_{p\mid d}\frac{p-|A_p|}{|A_p|}\notag\\
 &=|A|^2\prod_{p\mid Q}\left(|A_p|-1+\frac{p-|A_p|}{|A_p|}\right)
 \notag\\
 &=|A|^2\prod_{p\mid Q}\frac{(|A_p|-1)^2+p-1}{|A_p|}.\label{eq:weighted-first}
\end{align}

We next compute $\sum_t w(t)^2$.
By the Chinese remainder theorem, as $t$ ranges over the residues
modulo $Q$, the tuple $(t\bmod p)_{p\mid Q}$ ranges over all choices
of residues modulo the primes dividing $Q$, each exactly once.
Since $w_p(t)$ depends only on $t\bmod p$, the sum factors:
\begin{align}
 \sum_{t\bmod Q}w(t)^2
 &=\prod_{p\mid Q}\left(\sum_{u\bmod p}w_p(u)^2\right)\notag\\
 &=\prod_{p\mid Q}\bigl((p-1)^2+(p-1)(|A_p|-1)^2\bigr)\notag\\
 &=\prod_{p\mid Q}(p-1)\bigl((|A_p|-1)^2+p-1\bigr).
 \label{eq:weighted-second}
\end{align}
Here the second equality uses $w_p(0)=p-1$ and
$w_p(u)=|A_p|-1$ for the $p-1$ nonzero residues $u\bmod p$.

Substituting~\eqref{eq:weighted-first} and~\eqref{eq:weighted-second}
into~\eqref{eq:weighted-cs} gives
\begin{align*}
 E_Q(A)
 &\ge |A|^4\prod_{p\mid Q}\frac{(|A_p|-1)^2+p-1}{(p-1)|A_p|^2}\\
 &=\frac{|A|^4}{Q}\prod_{p\mid Q}
   \left(1+\frac{(p/|A_p|-1)^2}{p-1}\right).
\end{align*}
This proves the proposition.
\end{proof}

\section{Upper bounds on $E_{Q}(A)$}

A nonempty $B_k$-set is a $B_j$-set for every $j\le k$: one can pad
an equality of $j$-term sums by the same $k-j$ terms on both sides.
In particular, the sets considered below are Sidon.

\begin{lemma}\label{lem:upper}
Let $A\subseteq[N]$, and let $Q$ be a positive integer.
If $A$ is a $B_3$-set, then
\begin{equation}\label{eq:upper-three}
 E_Q(A)\le 2|A|^2+\frac{16N|A|}{Q}.
\end{equation}
If $A$ is a $B_4$-set, then
\begin{equation}\label{eq:upper-four}
 E_Q(A)\le 2|A|^2+4|A|^3+\frac{16N}{Q}.
\end{equation}
\end{lemma}

\begin{proof}
For an integer $t$, let
\[
 f(t)=\#\{(a,b,c,d)\in A^4:a+b-c-d=t\}.
\]
Relabeling the four variables shows that
\begin{equation}\label{eq:energy-integer}
 E_Q(A)=\sum_{\substack{|t|\le2N\\Q\mid t}}f(t).
\end{equation}
We first bound $f(t)$ in three cases.

\smallskip
\noindent\emph{The case $t=0$.}
Since $A$ is Sidon, the positive and negative summands must agree
as multisets. Hence $f(0)=2|A|^2-|A|\le2|A|^2$.

\smallskip
\noindent\emph{The case $t\in(A-A)\setminus\{0\}$.}
The Sidon property gives a unique representation $t=x-y$ with
$x,y\in A$. Any representation $a+b-c-d=x-y$ gives
\[
 a+b+y=c+d+x.
\]
The $B_3$ property forces equality of the multisets
$\{a,b,y\}$ and $\{c,d,x\}$. Because $x\ne y$, one of $a,b$
equals $x$, one of $c,d$ equals $y$, and the remaining two entries
are equal. There are at most $4|A|$ possibilities, so $f(t)\le4|A|$.

\smallskip
\noindent\emph{The case $t\notin A-A$.}
Fix $a\in A$. If
\[
 a+b-c-d=a+b'-c'-d'=t,
\]
then $b+c'+d'=b'+c+d$. The $B_3$ property gives equality of
these multisets. Since $b\notin\{c,d\}$ (otherwise $t\in A-A$),
we must have $b=b'$, and then $\{c,d\}=\{c',d'\}$.
There are at most two ordered triples $(b,c,d)$ for each $a$,
so $f(t)\le2|A|$.

If $A$ is a $B_4$-set, this last estimate improves to $f(t)\le4$.
Indeed, two representations of $t$ imply
\[
 \{a,b,c',d'\}=\{a',b',c,d\}
 \quad\text{as multisets}.
\]
Neither representation has a common entry among its positive and
negative summands, since $t\notin A-A$. Thus the multiset $\{a,b\}$
must equal $\{a',b'\}$, and $\{c,d\}$ must equal $\{c',d'\}$.
Only the at most four orderings remain.

There are at most $4N/Q$ nonzero multiples of $Q$ in $[-2N,2N]$.
For a $B_3$-set, every such $t$ has $f(t)\le4|A|$, proving
\eqref{eq:upper-three}. For a $B_4$-set, the nonzero elements of
$A-A$ contribute at most $4|A|\cdot |A|^2=4|A|^3$, and all remaining
nonzero multiples contribute at most $4\cdot4N/Q$. Together with
the $t=0$ term, this proves~\eqref{eq:upper-four}.
\end{proof}

\begin{proof}[Proof of \Cref{thm:main}]
We may assume $|A|\ge1$. Let $S$ be the subset of the $s$ smallest primes of $\cP$, and set
\[
 Q=\prod_{p\in S}p,
 \qquad
 L(S)=\prod_{p\in S}
 \left(1+\frac{(\alpha^{-1}-1)^2}{p-1}\right).
\]
By \Cref{prop:lower}, $E_Q(A)\ge |A|^4L(S)/Q$.
Combining this with \Cref{lem:upper} gives
\begin{align}
 L(S)\frac{|A|^3}{N}
 &\le16+\frac{2Q|A|}{N} &&(k=3),\label{eq:comparison-three}\\
 L(S)\frac{|A|^4}{N}
 &\le16+\frac{2Q|A|^2+4Q|A|^3}{N} &&(k=4).\label{eq:comparison-four}
\end{align}
Keep $s$, and hence $S$ and $Q$, fixed. The elementary estimate
$|A|=O_k(N^{1/k})$ from~\eqref{eq:counting} yields
\[
\frac{|A|}{N^{1/k}}
\le \left(\frac{16+o_s(1)}{L(S)}\right)^{1/k}.
\]
Finally, $L(S)\to\infty$ as $s\to\infty$, since
the logarithm of its $p$-factor is
$(\alpha^{-1}-1)^2/p+O_\alpha(p^{-2})$ and
$\sum_{p\in\cP}1/p=\infty$.
Thus by choosing $s$ arbitrarily large and letting $N\to \infty$, we obtain $|A|=o(N^{1/k})$.
\end{proof}

\begin{proof}[Proof of Corollary~\ref{cor:powers}]
We use the following classical fact.

\begin{lemma}[{\cite[Corollary~4.12(c)]{MV}}]\label{lem:primes}
The primes congruent to $1$ modulo $3$ satisfy
\[
 \sum_{p\equiv1\; (\mathrm{mod}\;3)}\frac1p=\infty.
\]
\end{lemma}

Suppose first that $k=3$. For every prime $p\equiv1\pmod3$, the
map $x\mapsto x^3$ on the cyclic group $\F_p^\times$ has a kernel
of size $3$, so its image has size $(p-1)/3$. Including zero gives $|A_p|\le1+\frac{p-1}{3}=\frac{p+2}{3}\le\frac p2$. The last inequality holds since these primes are at least $7$.
Apply \Cref{thm:main} with $N=n^3$, $\alpha=1/2$, and
$\cP=\{p:p\equiv1\pmod3\}$, using Lemma~\ref{lem:primes}.

Now suppose that $k=4$. Every fourth power is a square, so for
every odd prime $p$, $|A_p|\le\frac{p+1}{2}\le\frac{2p}{3}$. The reciprocal sum over the odd primes diverges (already by
Lemma~\ref{lem:primes}). Apply \Cref{thm:main} with $N=n^4$,
$\alpha=2/3$, and $\cP$ the odd primes. Both applications give
$|A|=o(n)$.
\end{proof}

\section{Concluding remarks}

For $k\ge5$, the $B_4$ upper bound still applies, but the comparison
controls only the fourth power of $|A|$:
\[
 |A|^4L(S)\le16N+2Q|A|^2+4Q|A|^3.
\]
Thus the argument gives only $|A|=o(N^{1/4})$, weaker than the
elementary bound $|A|=O_k(N^{1/k})$.

A natural attempt to exploit the stronger $B_k$ hypothesis is to
replace $E_Q(A)=E_{2,Q}(A)$ by a higher energy $E_{r,Q}(A)$,
which counts pairs of ordered $r$-term sums congruent modulo $Q$.
For $r\ge2$ and squarefree $Q$ supported on $\cP$, it is not difficult to generalize Proposition \ref{prop:lower} to show that
\[
 E_{r,Q}(A)\ge\frac{|A|^{2r}}Q
 \prod_{p\mid Q}
 \left(1+\frac{(\alpha^{-1}-1)^r}{(p-1)^{r-1}}\right).
\]
For $r=2$, recall that the product becomes arbitrarily large because
$\sum_{p\in\cP}1/p=\infty$. For each fixed $r\ge3$, however,
it stays bounded uniformly in $Q$, even when all primes are allowed,
since each factor is $1+O_{\alpha,r}(p^{1-r})$.
Thus this estimate gives only a bounded gain over $|A|^{2r}/Q$,
whereas our little-$o$ conclusion relies on making that gain
arbitrarily large.

Despite this obstruction, we expect ill-distribution to force a saving for every $k$.

\begin{conjecture}\label{conj:ill-distributed}
Fix an integer $k\ge5$ and a real number $\alpha$ such that $1/k\le\alpha<1$. If $A\subseteq[N]$ is a $B_k$-set
satisfying $|A_p|\le\alpha p$ for every prime $p$, then
\[
 |A|=o_{k,\alpha}(N^{1/k})
 \qquad(N\to\infty).
\]
\end{conjecture}

The analogous assertion holds for $k=2,3,4$ by
\eqref{Sidon} and Theorem~\ref{thm:main}. The range $\alpha<1/k$ follows directly from Gallagher's
larger sieve. At the endpoint $\alpha=1/k$, the larger-sieve
bound and the elementary $B_k$ bound coincide, and Conjecture \ref{conj:ill-distributed} predicts that their simultaneous hypotheses force a little $o$ improvement.

For sets of perfect powers, in forthcoming work we use additional arithmetic structure to prove the corresponding bound
for every fixed integer $k\ge2$: every $B_k$-set
$A\subseteq\{1^k,2^k,\ldots,n^k\}$ satisfies $|A|=o(n)$ as
$n\to\infty$.

\section{AI declaration}
The final proof of Theorem \ref{thm:main} was discovered through repeated interactions with Chat GPT-6 Astra. We would also like to acknowledge the usage of AI in editing the manuscript. The authors are entirely responsible for its content.

\end{document}